\documentclass{amsart}
\usepackage{amsfonts,amssymb,amsmath,amsthm}
\usepackage{url}
\usepackage{enumerate}

\usepackage{graphics}

\newtheorem{thrm}{Theorem}[section]
\newtheorem{lem}[thrm]{Lemma}

\theoremstyle{definition}

\numberwithin{equation}{section}

\author{Adrian \L ydka }
\address{Institute of Mathematics and Cryptology\\
Military University of Technology  \\
Kaliskiego 2, 00-908 Warsaw, Poland}
\email{adrian.lydka@wat.edu.pl}

\begin{document}

\title[On some properties of the function ...]{On some properties of the function of the number of relatively prime subsets of $\{1, 2, ..., n\}$}

\begin{abstract}
In the paper we solve few problems proposed by Prapanpong Pongsriiam. 
Let $f(n)$ denote the number of relatively prime subsets of $\{1, 2, 3, \dots, n\}$ and $g(n)$ denote the number of subsets $A$ of $\{1, 2, 3, \dots, n\}$ such that gcd$(A)>1$ and gcd$(A, n+1)=1$ . We show that $f_n^2-f_{n-k}f_{n+k}>0$ for $n\geq k+1\quad (k\geq2)$. We also show  $\frac{g(6n-2)}{g(6n-4)}>\frac{g(6n)}{g(6n-2)}>\frac{g(6n+2)}{g(6n)}<\frac{g(6n+4)}{g(6n+2)}$ for large $n$.


\end{abstract}

\maketitle

\section{Introduction}

A finite set $A$ is said to be relatively prime if $\text{gcd} (A) = 1$.

Let $f X)$ denote the number of relatively prime subsets of $X$.

Let $f(n)$ be the number of relatively prime subsets of $\{1, 2, 3, \dots, n\}$ in other words $f(n)=f([1, n])$. Sometimes we write $f_n$ instead of $f(n)$.

Moreover, define function $g(n)$ by formula

\begin{equation}\label{funkcjag}
g(n)=\sum_{\substack{\emptyset\neq A\subseteq [1, n] \\\text{gcd}(A)>1\\\text{gcd}(A, n+1)=1}}1.
\end{equation}

We will use two inequalities

\begin{lem}[\cite{nathanson1}, Theorem 2]\label{nath}
\begin{equation}\label{gol}
2^n-2^{\left\lfloor\frac{n}{2}\right\rfloor}-n2^{\left\lfloor\frac{n}{3}\right\rfloor}\leq f(n)\leq 2^n-2^{\left\lfloor\frac{n}{2}\right\rfloor}.
\end{equation}
\end{lem}

Moreover, we know that (Lemma 4 in \cite{pongsriiam1})

\begin{equation}\label{gol2}
g(n)=\sum_{\substack{2\leq d\leq n \\ (d, n+1)=1}}f\left(\left\lfloor\frac{n}{d}\right\rfloor\right).
\end{equation}

and 

\begin{equation}\label{fwzor}
f(n)=\sum_{d\leq n}\mu(d)\left(2^{\left\lfloor\frac{n}{d}\right\rfloor}-1\right).
\end{equation}

More information on the function $f_n$ can be found in the sequence A085945 in \cite{sloane}.

In paper \cite{pongsriiam1} Pongsriiam proved that $f_n^2-f_{n-1}f_{n+1}$ is positive for every odd number $n\geq3$ and negative for every even number $n$. 

Recall that a sequence $(a_n)_{n\geq0}$ is said to be log-concave if $a_n^2-a_{n-1}a_{n+1}>0$
for every $n>1$ and is said to be log-convex if $a_n^2-a_{n-1}a_{n+1}<0$ for every $n>1$.
Stirling numbers, Bessel numbers are examples of log-concave sequences. Some sequences are not log-concave, but have similar properties. For example, if $(F_n)_{n\geq0}$ is the Fibonacci sequence or $F_n=f_n$, then $F_n^2-F_n F_{n+1}= (-1)^{n-1}$ , which is positive for odd n and negative for even n (so called alternating sequence). In addition, the sequence $(f_n)_{n\geq1}$ seems to have strong log-property ((Recall that $(a_n)_{n\geq0}$ is said to be strong log-concave if $a_n^2-a_{n-k}a_{n+k}>0$ for
every $k\geq1$ and $n>k$)). For example, in the paper \cite{pongsriiam1} Pongsriiam checked that $f_n^2-f_{n-2}f_{n+2}>0 (\text{for} 2<n\leq 50)$, $f_n^2-f_{n-3}f_{n+3}>0 (\text{for} 3<n\leq 50)$ and $f_n^2-f_{n-4}f_{n+4}>0 (\text{for} 4<n\leq 50)$. In our paper we prove that these inequalities are true for all $n>2, 3, 4,$ respectively.

In this paper we prove that $f_n^2-f_{n-k}f_{n+k}>0  (\text{for large } n\geq k+1 \text{ and } k\geq2)$

We also propose new term :almost strong log-concave sequence if $a_n^2-a_{n-k}a_{n+k}>0$ for
every $k\geq k_0$ and $n>k$) for some constant $k_0\geq2$.

In paper \cite{pongsriiam1} Pongsriiam also asked is it true that
\begin{equation}
\frac{g(6n-2)}{g(6n-4)}>\frac{g(6n)}{g(6n-2)}>\frac{g(6n+2)}{g(6n)}<\frac{g(6n+4)}{g(6n+2)}.
\end{equation}

In Section.\ref{gineq}. we prove above inequalities for large $n$.


\section{Sign of $f_n^2-f_{n-k}f_{n+k}>0$ for $n>k$ in general}

First, using formula  (\ref{fwzor}), we can write the following GP/PARI code :

$a(n)=sum(k=1, n, \text{moebius}(k)*(2^{\text{floor}}(n/k)-1))$

$for(n=6,50,print(a(n)^2-a(n-l)*a(n+l)))$.

$for (l=2..8)$

We obtain that inequality $f_n^2-f_{n-k}f_{n+k}>0$ is true for $k=2, 3, ..., 8$ and $n\leq 50$.

Using estimation (\ref{gol}) we get

\begin{equation}
\begin{split}
f_n^2&\geq\left(2^n-2^{\left\lfloor\frac{n}{2}\right\rfloor}-n2^{\left\lfloor\frac{n}{3}\right\rfloor}\right)^2\\
&=2^{2n}-2^{n+\left\lfloor\frac{n}{2}\right\rfloor+1}-n2^{n+\left\lfloor\frac{n}{3}\right\rfloor+1}+2^{2\cdot\left\lfloor\frac{n}{2}\right\rfloor}+n2^{\left\lfloor\frac{n}{2}\right\rfloor+\left\lfloor\frac{n}{3}\right\rfloor+1}+n^2 2^{2\cdot\left\lfloor\frac{n}{3}\right\rfloor}\\
\end{split}
\end{equation}

and

\begin{equation}
\begin{split}
f_{n-k}f_{n+k}&\leq\left(2^{n-k}-2^{\left\lfloor\frac{n-k}{2}\right\rfloor}\right)\left(2^{n+k}-2^{\left\lfloor\frac{n+k}{2}\right\rfloor}\right)\\
&=2^{2n}-2^{n+\left\lfloor\frac{n+k}{2}\right\rfloor-k}-2^{n+\left\lfloor\frac{n-k}{2}\right\rfloor+k}+2^{\left\lfloor\frac{n-k}{2}\right\rfloor+\left\lfloor\frac{n+k}{2}\right\rfloor}
\end{split}
\end{equation}

So

\begin{equation}
\begin{split}
f_n^2-f_{n-k}f_{n+k}&\geq 2^{n+\left\lfloor\frac{n+k}{2}\right\rfloor}-2^{n+\left\lfloor\frac{n}{2}\right\rfloor+1}+2^{n-k+\left\lfloor\frac{n+k}{2}\right\rfloor}-n2^{n+\left\lfloor\frac{n}{3}\right\rfloor+1}\\
&\quad+2^{2\cdot\left\lfloor\frac{n}{2}\right\rfloor}-2^{\left\lfloor\frac{n+k}{2}\right\rfloor+\left\lfloor\frac{n-k}{2}\right\rfloor}+n2^{\left\lfloor\frac{n}{2}\right\rfloor+\left\lfloor\frac{n}{3}\right\rfloor+1}+n^2 2^{2\cdot\left\lfloor\frac{n}{3}\right\rfloor}
\end{split}
\end{equation}

\subsection{Case k=2}

\begin{equation}
\begin{split}
f_{n-2}f_{n+2}&\leq\left(2^{n-2}-2^{\left\lfloor\frac{n-2}{2}\right\rfloor}\right)\left(2^{n+2}-2^{\left\lfloor\frac{n+2}{2}\right\rfloor}\right)\\
&=2^{2n}-2^{n+\left\lfloor\frac{n}{2}\right\rfloor+1}-2^{n+\left\lfloor\frac{n}{2}\right\rfloor-1}+2^{2\cdot\left\lfloor\frac{n}{2}\right\rfloor}
\end{split}
\end{equation}

We show that $f_n^2-f_{n-2}f_{n+2}$ for $n\geq51$.

\begin{equation}
\begin{split}
f_n^2-f_{n-2}f_{n+2}   &\geq2^{n+\left\lfloor\frac{n}{2}\right\rfloor-1}-n2^{n+\left\lfloor\frac{n}{3}\right\rfloor+1}+n2^{\left\lfloor\frac{n}{2}\right\rfloor+\left\lfloor\frac{n}{3}\right\rfloor+1}+n^2 2^{2\cdot\left\lfloor\frac{n}{3}\right\rfloor}\\
&>2^{n+\left\lfloor\frac{n}{2}\right\rfloor-1}-n2^{n+\left\lfloor\frac{n}{3}\right\rfloor+1}>2^{n+\left\lfloor\frac{n}{2}\right\rfloor-1}-2^{n+\left\lfloor\frac{n}{3}\right\rfloor+\log_2 n+1}>0,
\end{split}
\end{equation}

because $\left\lfloor\frac{n}{2}\right\rfloor-\left\lfloor\frac{n}{3}\right\rfloor-\log_2 n-2\geq \frac{n-1}{2}-\frac{n}{3}-\log_2 n-2=\frac{n-6\log_2 n-15}{6}>0$ for $n\geq 51.$ (Consider function $h(x)=x-6\log_2 x-15$, $h(51)=36-6\log_2 55>0$, $h'(x)=1-\frac{6}{x\ln{2}}>0$ for $x\geq51\geq 6\ln{2}$).

\subsection{Case k=3}

\begin{equation}
\begin{split}
f_n^2-f_{n-3}f_{n+3}&\geq 2^{n+\left\lfloor\frac{n+1}{2}\right\rfloor+1}-2^{n+\left\lfloor\frac{n}{2}\right\rfloor+1}+2^{n-2+\left\lfloor\frac{n+1}{2}\right\rfloor}-n2^{n+\left\lfloor\frac{n}{3}\right\rfloor+1}\\
&\quad+2^{2\cdot\left\lfloor\frac{n}{2}\right\rfloor}-2^{\left\lfloor\frac{n+1}{2}\right\rfloor+\left\lfloor\frac{n-1}{2}\right\rfloor}+n2^{\left\lfloor\frac{n}{2}\right\rfloor+\left\lfloor\frac{n}{3}\right\rfloor+1}+n^2 2^{2\cdot\left\lfloor\frac{n}{3}\right\rfloor}
\end{split}
\end{equation}

We show that $f_n^2-f_{n-3}f_{n+3}$ for $n\geq51$.

First we see that

\begin{equation}
\begin{split}
2^{n+\left\lfloor\frac{n+1}{2}\right\rfloor+1}-2^{n+\left\lfloor\frac{n}{2}\right\rfloor+1}+2^{2\cdot\left\lfloor\frac{n}{2}\right\rfloor}-2^{\left\lfloor\frac{n+1}{2}\right\rfloor+\left\lfloor\frac{n-1}{2}\right\rfloor}>0
\end{split}
\end{equation}

So, it is enough to prove that 
\begin{equation}\label{dlatrzy}
2^{n-2+\left\lfloor\frac{n+1}{2}\right\rfloor}-n2^{n+\left\lfloor\frac{n}{3}\right\rfloor+1}\geq0
\end{equation}

for $n\geq51$.

We have $\left\lfloor\frac{n+1}{2}\right\rfloor-\left\lfloor\frac{n}{3}\right\rfloor-3-\log_2{n}\geq \frac{n}{2}-\frac{n}{3}-3-\log_2{n}$.

Consider function $h(x)=\frac{x}{2}-\frac{x}{3}-3-\log_2{x}$. We  have $h(51)=\frac{33}{6}-\log_2{51}>0$ and $h'(x)=\frac{1}{6}-\frac{1}{x\ln{2}}, h'(x)>0$ for $x>51>\frac{6}{\ln{2}}$.  So inequality (\ref{dlatrzy}) is true.

\subsection{Case k=4}

\begin{equation}
\begin{split}
f_n^2-f_{n-4}f_{n+4}&\geq 2^{n+\left\lfloor\frac{n}{2}\right\rfloor+2}-2^{n+\left\lfloor\frac{n}{2}\right\rfloor+1}+2^{n-2+\left\lfloor\frac{n}{2}\right\rfloor}-n2^{n+\left\lfloor\frac{n}{3}\right\rfloor+1}\\
&+n2^{\left\lfloor\frac{n}{2}\right\rfloor+\left\lfloor\frac{n}{3}\right\rfloor+1}+n^2 2^{2\cdot\left\lfloor\frac{n}{3}\right\rfloor}
\end{split}
\end{equation}

\begin{equation}
\begin{split}
f_n^2-f_{n-4}f_{n+4}&\geq 2^{n+\left\lfloor\frac{n}{2}\right\rfloor+1}+2^{n-2+\left\lfloor\frac{n}{2}\right\rfloor}-n2^{n+\left\lfloor\frac{n}{3}\right\rfloor+1}\\
&+n2^{\left\lfloor\frac{n}{2}\right\rfloor+\left\lfloor\frac{n}{3}\right\rfloor+1}+n^2 2^{2\cdot\left\lfloor\frac{n}{3}\right\rfloor}
\end{split}
\end{equation}

Now inequality (\ref{dlatrzy}) implies statement.

\subsection{Case k=5}

\begin{equation}
\begin{split}
f_n^2-f_{n-5}f_{n+5}&\geq 2^{n+\left\lfloor\frac{n+5}{2}\right\rfloor}-2^{n+\left\lfloor\frac{n}{2}\right\rfloor+1}+2^{n-5+\left\lfloor\frac{n+5}{2}\right\rfloor}-n2^{n+\left\lfloor\frac{n}{3}\right\rfloor+1}\\
&\quad+2^{2\cdot\left\lfloor\frac{n}{2}\right\rfloor}-2^{\left\lfloor\frac{n+5}{2}\right\rfloor+\left\lfloor\frac{n-5}{2}\right\rfloor}+n2^{\left\lfloor\frac{n}{2}\right\rfloor+\left\lfloor\frac{n}{3}\right\rfloor+1}+n^2 2^{2\cdot\left\lfloor\frac{n}{3}\right\rfloor}
\end{split}
\end{equation}

We prove that for $n\geq 36$ above term is positive.

This term is great than

\begin{equation}
2^{n+\left\lfloor\frac{n+1}{2}\right\rfloor+1}-n2^{n+\left\lfloor\frac{n}{3}\right\rfloor+1}=2^{n+\left\lfloor\frac{n}{3}\right\rfloor+1}\left(2^{\left\lfloor\frac{n+1}{2}\right\rfloor-\left\lfloor\frac{n}{3}\right\rfloor}-2^{\log_2{n}}\right)
\end{equation}

$\left\lfloor\frac{n+1}{2}\right\rfloor-\left\lfloor\frac{n}{3}\right\rfloor-\log_2{n}\geq \frac{n}{2}-\frac{n}{3}-\log_2{n}\geq 0$, for $n\geq 36$.

\subsection{Case k=6}

\begin{equation}
\begin{split}
f_n^2-f_{n-6}f_{n+6}&\geq 2^{n+\left\lfloor\frac{n+6}{2}\right\rfloor}-2^{n+\left\lfloor\frac{n}{2}\right\rfloor+1}+2^{n-6+\left\lfloor\frac{n+6}{2}\right\rfloor}-n2^{n+\left\lfloor\frac{n}{3}\right\rfloor+1}\\
&\quad+2^{2\cdot\left\lfloor\frac{n}{2}\right\rfloor}-2^{\left\lfloor\frac{n+6}{2}\right\rfloor+\left\lfloor\frac{n-6}{2}\right\rfloor}+n2^{\left\lfloor\frac{n}{2}\right\rfloor+\left\lfloor\frac{n}{3}\right\rfloor+1}+n^2 2^{2\cdot\left\lfloor\frac{n}{3}\right\rfloor}
\end{split}
\end{equation}

We prove that for $n\geq 36$ above term is positive.

This term is great than

\begin{equation}
2^{n+\left\lfloor\frac{n}{2}\right\rfloor+2}-n2^{n+\left\lfloor\frac{n}{3}\right\rfloor+1}=2^{n+\left\lfloor\frac{n}{3}\right\rfloor+1}\left(2^{\left\lfloor\frac{n}{2}\right\rfloor-\left\lfloor\frac{n}{3}\right\rfloor+1}-2^{\log_2{n}}\right),
\end{equation}

but $\left\lfloor\frac{n}{2}\right\rfloor-\left\lfloor\frac{n}{3}\right\rfloor+1-\log_2{n}\geq \frac{n-1}{2}-\frac{n}{3}+1-\log_2{n}\geq 0$, for $n\geq 36$.

\subsection{Case k=7}

\begin{equation}
\begin{split}
f_n^2-f_{n-7}f_{n+7}&\geq 2^{n+\left\lfloor\frac{n+7}{2}\right\rfloor}-2^{n+\left\lfloor\frac{n}{2}\right\rfloor+1}+2^{n-7+\left\lfloor\frac{n+7}{2}\right\rfloor}-n2^{n+\left\lfloor\frac{n}{3}\right\rfloor+1}\\
&\quad+2^{2\cdot\left\lfloor\frac{n}{2}\right\rfloor}-2^{\left\lfloor\frac{n+7}{2}\right\rfloor+\left\lfloor\frac{n-7}{2}\right\rfloor}+n2^{\left\lfloor\frac{n}{2}\right\rfloor+\left\lfloor\frac{n}{3}\right\rfloor+1}+n^2 2^{2\cdot\left\lfloor\frac{n}{3}\right\rfloor}
\end{split}
\end{equation}

We prove that for $n\geq 36$ above term is positive.

This term is great than

\begin{equation}
2^{n+\left\lfloor\frac{n+1}{2}\right\rfloor+2}-n2^{n+\left\lfloor\frac{n}{3}\right\rfloor+1}=2^{n+\left\lfloor\frac{n}{3}\right\rfloor+1}\left(2^{\left\lfloor\frac{n+1}{2}\right\rfloor-\left\lfloor\frac{n}{3}\right\rfloor+1}-2^{\log_2{n}}\right),
\end{equation}

but $\left\lfloor\frac{n+1}{2}\right\rfloor-\left\lfloor\frac{n}{3}\right\rfloor+1-\log_2{n}\geq \frac{n}{2}-\frac{n}{3}+1-\log_2{n}\geq 0$, for $n\geq 36$.

\subsection{Case k=8}

\begin{equation}
\begin{split}
f_n^2-f_{n-8}f_{n+8}&\geq 2^{n+\left\lfloor\frac{n+8}{2}\right\rfloor}-2^{n+\left\lfloor\frac{n}{2}\right\rfloor+1}+2^{n-8+\left\lfloor\frac{n+8}{2}\right\rfloor}-n2^{n+\left\lfloor\frac{n}{3}\right\rfloor+1}\\
&\quad+2^{2\cdot\left\lfloor\frac{n}{2}\right\rfloor}-2^{\left\lfloor\frac{n+8}{2}\right\rfloor+\left\lfloor\frac{n-8}{2}\right\rfloor}+n2^{\left\lfloor\frac{n}{2}\right\rfloor+\left\lfloor\frac{n}{3}\right\rfloor+1}+n^2 2^{2\cdot\left\lfloor\frac{n}{3}\right\rfloor}
\end{split}
\end{equation}

We prove that for $n\geq 36$ above term is positive.

This term is great than

\begin{equation}
2^{n+\left\lfloor\frac{n}{2}\right\rfloor+3}-n2^{n+\left\lfloor\frac{n}{3}\right\rfloor+1}=2^{n+\left\lfloor\frac{n}{3}\right\rfloor+1}\left(2^{\left\lfloor\frac{n}{2}\right\rfloor-\left\lfloor\frac{n}{3}\right\rfloor+2}-2^{\log_2{n}}\right),
\end{equation}

but $\left\lfloor\frac{n}{2}\right\rfloor-\left\lfloor\frac{n}{3}\right\rfloor+2-\log_2{n}\geq \frac{n-1}{2}-\frac{n}{3}+2-\log_2{n}\geq 0$, for $n\geq 36$.

\subsection{Case $k\geq 9$}

\begin{equation}
\begin{split}
f_n^2-f_{n-k}f_{n+k}&\geq 2^{n+\left\lfloor\frac{n+k}{2}\right\rfloor}-2^{n+\left\lfloor\frac{n}{2}\right\rfloor+1}+2^{n-k+\left\lfloor\frac{n+k}{2}\right\rfloor}-n2^{n+\left\lfloor\frac{n}{3}\right\rfloor+1}\\
&\quad+2^{2\cdot\left\lfloor\frac{n}{2}\right\rfloor}-2^{\left\lfloor\frac{n+k}{2}\right\rfloor+\left\lfloor\frac{n-k}{2}\right\rfloor}+n2^{\left\lfloor\frac{n}{2}\right\rfloor+\left\lfloor\frac{n}{3}\right\rfloor+1}+n^2 2^{2\cdot\left\lfloor\frac{n}{3}\right\rfloor}
\end{split}
\end{equation}

We prove that for $n\geq k+1$ above term is positive.

This term is great than

\begin{equation}
\begin{split}
2^{n+\left\lfloor\frac{n+k}{2}\right\rfloor-1}-n2^{n+\left\lfloor\frac{n}{3}\right\rfloor+1}
\end{split}
\end{equation}

It is enough to show that 

\begin{equation}
\begin{split}
(n+\left\lfloor\frac{n+k}{2}\right\rfloor-1)-(n+\left\lfloor\frac{n}{3}\right\rfloor+1+\log_2{n})\geq0,
\end{split}
\end{equation}

but

\begin{equation}
\begin{split}
(n+\left\lfloor\frac{n+k}{2}\right\rfloor-1)-(n+\left\lfloor\frac{n}{3}\right\rfloor+1+\log_2{n})\geq \frac{n+k-1}{2}-1-\frac{n}{3}-1-\log_2{n}\\
=\frac{n}{6}+\frac{k-5}{2}-\log_2{n}
\end{split}
\end{equation}

Let $i(x)=\frac{x}{6}+\frac{k-5}{2}-\log_2{x}$, $i(k+1)=\frac{4k-14}{6}-\log_2(k+1)$. Function i(x) is increasing for $x\geq8$ and $i(k+1)>0$ for $k\geq 9$.

\section{Proof that $\frac{g(6n-2)}{g(6n-4)}>\frac{g(6n)}{g(6n-2)}>\frac{g(6n+2)}{g(6n)}<\frac{g(6n+4)}{g(6n+2)}$ for large $n$.}\label{gineq}

\subsection{$\frac{g(6n-2)}{g(6n-4)}>\frac{g(6n)}{g(6n-2)}$}

Above inequality is equivalent to inequality

\begin{equation}\label{gol3}
[g(6n-2)]^2>g(6n)g(6n-4).
\end{equation}

Using (\ref{gol2}) we get estimation

\begin{equation}
\begin{split}
g(6n-2)&=\sum_{\substack{2\leq d\leq 6n-2 \\ (d, 6n-1)=1}}f\left(\left\lfloor\frac{6n-2}{d}\right\rfloor\right)=f(3n-1)+f(2n-1)\\&\quad+f\left(\left\lfloor\frac{3n-1}{2}\right\rfloor\right)+\chi_5(6n-1)f\left(\left\lfloor\frac{6n-2}{5}\right\rfloor\right)+C_1(n),
\end{split}
\end{equation}

where $C_1(n)=\sum\limits_{\substack{6\leq d\leq 6n-2 \\ (d, 6n-1)=1}}f\left(\left\lfloor\frac{6n-2}{d}\right\rfloor\right)\leq 6(n-1)f(n-1)\leq 6(n-1)2^{n-1}$.

Subsequently

\begin{equation}
\begin{split}
g(6n-4)&=\sum_{\substack{2\leq d\leq 6n-4 \\ (d, 6n-3)=1}}f\left(\left\lfloor\frac{6n-4}{d}\right\rfloor\right)=f(3n-2)+f\left(\left\lfloor\frac{3n}{2}\right\rfloor-1\right)\\&\quad+\chi_5(6n-3)f\left(\left\lfloor\frac{6n-4}{5}\right\rfloor\right)+C_2(n),
\end{split}
\end{equation}

where $C_2(n)=\sum\limits_{\substack{6\leq d\leq 6n-4 \\ (d, 6n-3)=1}}f\left(\left\lfloor\frac{6n-4}{d}\right\rfloor\right)\leq 6(n-1)f(n-1)\leq 6(n-1)2^{n-1}$.

Similarly

\begin{equation}
\begin{split}
g(6n)&=\sum_{\substack{2\leq d\leq 6n \\ (d, 6n+1)=1}}f\left(\left\lfloor\frac{6n}{d}\right\rfloor\right)=f(3n)+f(2n)+f\left(\left\lfloor\frac{3n}{2}\right\rfloor\right)\\&\quad+\chi_5(6n+1)f\left(\left\lfloor\frac{6n}{5}\right\rfloor\right)+f(n)+C_3(n),
\end{split}
\end{equation}

where $C_3(n)=\sum\limits_{\substack{7\leq d\leq 6n \\ (d, 6n+1)=1}}f\left(\left\lfloor\frac{6n}{d}\right\rfloor\right)\leq 6(n-1)f(n-1)\leq 6(n-1)2^{n-1}$.

Now, (\ref{gol3}) is equivalent to inequality

\begin{equation}
\begin{split}
&\left[f(3n-1)+f(2n-1)+f\left(\left\lfloor\frac{3n-1}{2}\right\rfloor\right)+\chi_5(6n-1)f\left(\left\lfloor\frac{6n-2}{5}\right\rfloor\right)+C_1(n)\right]^2\\
&>\left[f(3n)+f(2n)+f\left(\left\lfloor\frac{3n}{2}\right\rfloor\right)+\chi_5(6n+1)f\left(\left\lfloor\frac{6n}{5}\right\rfloor\right)+f(n)+C_3(n)\right]\\
&\quad\cdot\left[f(3n-2)+f\left(\left\lfloor\frac{3n}{2}\right\rfloor-1\right)+\chi_5(6n-3)f\left(\left\lfloor\frac{6n-4}{5}\right\rfloor\right)+C_2(n)\right]
\end{split}
\end{equation}

 Above inequality after calculation, cancellation summands which are $O\left(2^{\frac{9}{2}n}\right)$, leads to inequality

\begin{equation}
\begin{split}
[f(3n-1)]^2+2f(3n-1)f(2n-1)>f(3n)f(3n-2)+f(3n-2)f(2n),
\end{split}
\end{equation}
   
which is true for large $n$. So the inequality $\frac{g(6n-2)}{g(6n-4)}>\frac{g(6n)}{g(6n-2)}$ is true for large $n$.

\subsection{$\frac{g(6n)}{g(6n-2)}>\frac{g(6n+2)}{g(6n)}$}

Above inequality is equivalent to inequality

\begin{equation}\label{gol4}
[g(6n)]^2>g(6n-2)g(6n+2).
\end{equation}

 Above inequality after calculation, cancellation summands which are $O\left(2^{\frac{9}{2}n}\right)$, leads to inequality

\begin{equation}
\begin{split}
[f(3n)]^2+2f(3n)f(2n)>f(3n-1)f(3n+1)+f(2n-1)f(3n+1),
\end{split}
\end{equation}

which is true for large $n$.

\subsection{$\frac{g(6n+2)}{g(6n)}<\frac{g(6n+4)}{g(6n+2)}$}

Above inequality is equivalent to inequality

\begin{equation}\label{gol5}
[g(6n+2)]^2>g(6n)g(6n+4).
\end{equation}

 Above inequality after calculation, cancellation summands which are $O\left(2^{\frac{9}{2}n}\right)$, leads to inequality

\begin{equation}
\begin{split}
[f(3n+1)]^2<f(3n)f(2n+1)+f(2n)f(3n+2),
\end{split}
\end{equation}

Using inequalities from Lemma.\ref{nath}. we can prove that above inequality is true for large $n$.

We have proved that exists natural  number  $n_0$ such that inequality from the title of the section is true for $n\geq n_0$
. Exact value $n_0$ needs more careful calculation.



\end{document}